\documentclass[12pt]{article}
\usepackage{amssymb}
\usepackage{amsmath}
\usepackage{latexsym}
\usepackage{mathrsfs}
\usepackage{color}

\numberwithin{equation}{section}
\newtheorem{thm}{Theorem}
\newtheorem{prop}[thm]{Proposition}

\newtheorem{rem}[thm]{Remark}
\newtheorem{lem}[thm]{Lemma}
\newtheorem{corol}[thm]{Corollary}

\newcommand{\avint}{{\int\hspace{-12pt}-}}

\newcommand{\R}{{\mathbb{R}}}
\newcommand{\cC}{{\cal C}}

\newcommand{\bLip}{{\rm bLip}}
\newcommand{\Lip}{{\rm Lip}}

\newcommand{\darr}[4]{{\left\{\begin{array}{ll}
  {#1}&{#2}\\[0.2cm]
  {#3}&{#4}
\end{array}\right.}}

\def\Frac{\displaystyle\frac}

\title{Shape sensitivity analysis of the Hardy constant}

\author{
Gerassimos Barbatis\footnote{Department of Mathematics,
 University of Athens,  15784 Athens, Greece}\ \ 
and Pier Domenico Lamberti 
 \footnote{Dipartimento di Matematica,
Universit\`{a} degli Studi di Padova,
Via Trieste 63, 35121 Padova, Italy}}

\date{\ }

\begin{document}

\newcommand{\rea}{\mathbb{R}}

\maketitle

%
%
%

\noindent

%
%

\begin{abstract}
\noindent  We consider  the Hardy constant associated with a domain in the $n$-di\-men\-sio\-nal Euclidean space and we study its variation upon perturbation of the domain. We prove a Fr\'{e}chet differentiability result and establish a Hadamard-type formula for the corresponding derivatives. We also prove 
a stability result for the minimizers of the Hardy quotient. Finally, we prove stability estimates in terms of the Lebesgue measure of the symmetric 
difference of domains. 

\end{abstract}

\vspace{11pt}


\vspace{6pt}

\parskip=5pt

\section{Introduction}

Let $\Omega $ be a bounded domain in ${\mathbb{R}}^n$, $d_{\Omega }(x)= {\rm dist }(x, \partial \Omega)$, $x\in \Omega$, and 
$p\in ]1,\infty [$. If there exists $c>0$ such that 
\begin{equation}
\label{mainhardy}
\int_{\Omega }|\nabla u|^pdx \geq c \int_{\Omega }\frac{|u|^p}{d^p_{\Omega }}dx \, , \quad \mbox{ for all $u\in C^{\infty }_c(\Omega)$}, 
\end{equation}
we then say that the $L^p$ Hardy inequality holds in $\Omega$.
The best constant for inequality (\ref{mainhardy})
is called the $L^p$ Hardy constant of $\Omega$ and we shall denote it by $H_p(\Omega)$. It is well-known that if $\Omega $ is regular enough then the $L^p$ Hardy inequality is valid for all $p \in ]1,\infty[$; moreover if $\Omega$ is convex, and more generally if it is weakly mean convex, i.e. if $\Delta d_{\Omega}\le 0$ in the distributional sense in $\Omega$, then $H_p(\Omega )=((p-1)/p)^p$. 

The study of inequality (\ref{mainhardy}) has a long history which goes back to Hardy himself, see \cite{permalkuf}. In the last twenty years there has been a growing interest in the study of
Hardy inequalities, the existence and behavior of minimizers \cite{MMP, mash}, improved inequalities \cite{bremar, barfilter}, higher order analogues and other related problems. 

The precise evaluation of $H_p(\Omega )$ for domains $\Omega$ that are not weakly  mean convex is a difficult problem. 
There are only few examples of such domains for which $H_p(\Omega )$
is known and these are only for the case $p=2$ and for very special
domains $\Omega$. 
 Even the problem of estimating from below $H_p(\Omega)$ is difficult and most results again are for $p=2$. One such  result is the well known theorem by A. Ancona which states that $H_2(\Omega )\geq 1/16$ for all simply connected planar domains. We refer to  \cite{davies, MMP, barfilter, laso, BarTer, avk} for more information on the Hardy constant. 

In this paper we study the variation of $H_p(\Omega )$ upon variation of the domain $\Omega$. 
This probem can be considered as a spectral perturbation problem. Indeed, if there exists a minimizer $u\in W^{1,p}_0(\Omega )$  for the Hardy quotient associated with (\ref{mainhardy}) then $u$ is a solution to the equation
\begin{equation}\label{pde}
-\Delta_p u=H_p(\Omega )\frac{|u|^{p-2}u}{d_{\Omega }^p}
\end{equation}
where $\Delta _pu={\rm div}(|\nabla u|^{p-2}\nabla u)$ is the $p$-Laplacian.
Domain perturbation problems  have been extensively studied in the case of the Dirichlet Laplacian as well as for more general elliptic operators, such as operators satisfying other boundary conditions, higher order operators and operators with variable coefficients.
When studying such problems, there are broadly speaking two types of results: qualitative and quantitative. The former provide information such as continuity or analyticity, while the second involve stability properties, possibly together with related estimates. The relevant literature is vast, and we refer to \cite{bulamat, bulala, daners, hale, henry} and references therein for more information; in particular, for the  $p$-laplacian we refer to \cite{bulaplap, garsab, laplap}.

In this paper we obtain both qualitative and quantitative results on the domain dependence of $H_p(\Omega)$.
In  Theorem \ref{thm:diff}, we assume that $\Omega $ is of class $C^2$ with $H_p(\Omega )<((p-1)/p)^p$ and we establish the Fr\'{e}chet differentiability of $H_p(\phi (\Omega) )$ with respect to the $C^2$ diffeomorphism $\phi$. In particular we provide a Hadamard-type formula for the Fr\'{e}chet differential. For our proof we make essential use of certain results of \cite{MMP}, where it was shown in particular that if $H_p(\Omega )<((p-1)/p)^p$ then the Hardy
quotient admits a positive minimizer $u$ which behaves like $d^{\alpha}_{\Omega}$ near $\partial\Omega$ for a suitable $\alpha >0$.
In fact, in Theorem~\ref{stabmin} we also prove the stability of the minimizer $u$ in  $W^{1,p}_0(\Omega)$; this is of independent interest but is also used in the proof of Theorem \ref{thm:diff}.

We subsequently consider stability estimates for $H_p(\Omega)$. In Theorem \ref{thm:aaa} we prove under certain assumptions that the Hardy constant $H_p(\Omega)$ of a $C^2$ domain $\Omega$ is upper semicontinuous with respect to bi-Lipschitz tranformations $\phi$. In Theorem~\ref{thm:niki} we consider the stability of the Hardy constant when $\Omega$ is subject to a localized perturbation which transforms it to a  domain $\tilde \Omega$. Assuming that both $\Omega $ and $\tilde\Omega$ are of class $C^2$ we obtain stability estimates for the $L^p$ Hardy constant in terms of the Lebesgue measure of the symmetric difference $\Omega \triangle \tilde\Omega$. Estimates of this type have been recently obtained for eigenvalues of various classes of operators; we refer to \cite{BL1, bulala, bulamat} and references therein for more information.

The paper  is organized as follows. In Section 2 we introduce our notation and prove a   general Lipschitz continuity result. Section 3 is devoted to the proof of differentiability results, the Hadamard formula and the stability of minimizers. In Section 4 we prove stability estimates in terms of  the Lebesgue measure of
the symmetric difference of the domains.

\section{Preliminaries}

Let $\Omega$ be a bounded domain (i.e. a bounded connected open set) in $\R^n$. Given $p\in ]1,+\infty[$ we denote by $W^{1,p}_0(\Omega)$ the closure in the standard Sobolev space $W^{1,p}(\Omega)$ of the set of all smooth functions with compact support in $\Omega$.

If  $u\in W^{1,p}_0(\Omega)$, $u\ne 0$, we then denote by $R_{\Omega}[u]$ the Rayleigh quotient
\[
R_{\Omega,p}[u] =\frac{\int_{\Omega}|\nabla u|^pdx}{\int_{\Omega}\frac{|u|^p}{d_{\Omega}^p}dx},
\]
and we set
\begin{equation}
H_p(\Omega) =\inf_{u\in W^{1,p}_0(\Omega), u\neq 0} R_{\Omega,p}[u].
\label{ray}
\end{equation}

If $H_p(\Omega)>0$ we then say that the $L^p$ Hardy inequality is valid on
$\Omega$.  

It is well known that if $\Omega$ has a Lipschitz continuous boundary then $0<H_p(\Omega) \leq ((p-1)/p)^p$ and it has been proved in \cite{MMP,mash} that if $\Omega$ is of class $C^2$ then there exists a minimizer $u$ in (\ref{ray}) if and only if $H_p(\Omega)< ((p-1)/p)^p$; moreover,  such minimizer is unique up to a multiplicative constant, can be chosen to be positive and there exists $c>0$ such that
\begin{equation}\label{decay}
c^{-1}d_{\Omega}(x)^{\alpha} \leq u(x) \leq  c d_{\Omega}(x)^{\alpha},\ \ x\in \Omega ,
\end{equation}
where $\alpha>(p-1)/p$ is the largest solution to the equation 
\begin{equation}
(p-1)\alpha^{p-1}(1-\alpha)=H_p(\Omega).
\label{alpha}
\end{equation}

Given a Lipschitz map $\phi:\Omega\to\phi(\Omega)$ we define $\Lip(\phi)=\|\nabla\phi\|_{L^{\infty}(\Omega)}$. 
For $L>0$ we define the uniform class of bi-Lipschitz maps
\begin{eqnarray*}\lefteqn{
\bLip_L(\Omega)= \{ \phi:\Omega\to\phi(\Omega) : \; \phi,\ \phi^{(-1)}\  {\rm are\ Lipschitz\ continuous }}\\ & & \qquad\qquad\qquad\qquad\qquad\qquad\qquad\quad {\rm and}\ \Lip(\phi) ,\, \Lip(\phi^{(-1)}) \leq L\ \}.
\end{eqnarray*}

In the sequel we shall often use the fact that $H_p(\phi (\Omega ))$ depends continuously on $\phi $. In fact, we can prove the 
 following Lipschitz continuity result. 

Note that in the proof  of the following proposition as well as in the proofs of other statements in the sequel, by $c, c_1$ etc. we shall denote constants the value of which  may change
from line to line.

\begin{prop}
Let $\Omega$ be a bounded domain in $\R^n$, $p\in ]1,\infty [$ and $L>0$. There exists $c>0$ depending only on $n,p,L$ such that
\begin{equation}
| H_p(\phi(\Omega)) -  H_p(\Omega)|  \leq    c H_p(\Omega)\| \nabla\phi- I\|_{L^{\infty}(\Omega)} \, ,
\label{bbb}
\end{equation}
for all $\phi \in \bLip_L(\Omega)$ such that $\|\nabla\phi- I\|_{L^{\infty}(\Omega)}<c^{-1}$.
\label{thm:bbb}
\end{prop}
{\em Proof.} Let $u\in W^{1,p}_0(\Omega)$ be normalized by $\int_{\Omega}|u|^p/ d_{\Omega}^p \, dx=1$. Then $v:=u\circ \phi^{(-1)}$ belongs to $W^{1,p}_0(\phi(\Omega))$.
Changing variables we have
\begin{equation}
\frac{\int_{\phi(\Omega)}|\nabla v|^p dy}{\int_{\phi(\Omega)}\Frac{|v|^p}{d^p_{\phi(\Omega)}(y)}dy}
=\frac{\int_{\Omega}| (\nabla u) (\nabla\phi)^{-1}|^p \; |\det \nabla\phi | dx}{\int_{\Omega}\Frac{|u|^p}{d^p_{\phi(\Omega)}(\phi(x))}
|\det \nabla\phi | dx},
\label{chvar}
\end{equation}
so
\begin{eqnarray}
&& \frac{\int_{\phi(\Omega)}|\nabla v|^p dy}{\int_{\phi(\Omega)}\Frac{|v|^p}{d^p_{\phi(\Omega)}(y)}dy} -
\frac{\int_{\Omega}|\nabla u|^p dx}{\int_{\Omega}\Frac{|u|^p}{d^p_{\Omega}(x)}dx} \label{b4}\\
&=& \frac{ \int_{\Omega} \Big( |(\nabla u) (\nabla\phi)^{-1}|^p |\det \nabla\phi | -  |\nabla u|^p \Big) dx  -  
\int_{\Omega}|\nabla u|^p dx \; \Big(  \int_{\Omega}\Frac{|u|^p|\det \nabla\phi |}{d^p_{\phi(\Omega)}(\phi(x))}  dx \, - 1\Big)}
{\int_{\Omega}\Frac{|u|^p}{d^p_{\phi(\Omega)}(\phi(x))} |\det \nabla\phi | dx}.
\nonumber
\end{eqnarray}
Using the relation $||\det\nabla\phi| -1| \leq c|\nabla\phi -I|$, valid since $|\nabla\phi|< c_1$, we obtain after some simple computations that
\begin{equation}
\bigg|  | (\nabla u) (\nabla\phi)^{-1}|^p |\det \nabla\phi | -  |\nabla u|^p \bigg|  \leq c   |\nabla\phi -I| \cdot  |\nabla u|^p \, .   
\label{eq:bbc}
\end{equation}

 Note that $\phi $ admits a unique Lipschitz continuous extension on $\overline{\Omega}$.
Moreover, if $x,y\in {\mathbb{R}}^n$ are such that the ``open" line segment $]x,y[$ is
contained in $\Omega$, then
\begin{equation}
\Big|  |\phi(y)-\phi(x)| - |y-x| \Big|  \leq   \|\nabla\phi -I\|_{L^{\infty}(\Omega)} |y-x| \; .
\label{ster}
\end{equation}
Choosing $y$ to be a nearest boundary point to $x$ yields
\begin{equation}
   d_{\phi(\Omega)}(\phi(x)) \leq  d_{\Omega}(x)  \big( 1 +   \|\nabla\phi - I\|_{L^{\infty}(\Omega)}    \big) .
\label{b6}
\end{equation}
Hence, using also the normalization, we have 
\begin{equation}
\int_{\Omega} \frac{|u|^p}{d^p_{\phi(\Omega)}(\phi(x))} |\det \nabla\phi | dx  \geq  
 1  - c  \|\nabla\phi -I\|_{L^{\infty}(\Omega)}. 
\label{eq:bbc1}
\end{equation}
Combining (\ref{b4}), (\ref{eq:bbc}), (\ref{eq:bbc1}) and the normalization we conclude that
\begin{equation}
R_{\phi(\Omega)}[v]\leq (1+ c \|\nabla\phi - I\|_{L^{\infty}(\Omega)}) R_{\Omega}[u] ,
\label{ssoo}
\end{equation}
provided $ \|\nabla\phi - I\|_{L^{\infty}(\Omega)}$ is sufficiently small.
 Hence $H_p(\phi(\Omega)) \leq H_p(\Omega)(1 +c  \|\nabla\phi - I\|_{L^{\infty}(\Omega)})$.
Replacing $\Omega$ by $\phi(\Omega)$ and $\phi$ by $\phi^{(-1)}$ we obtain
$H_p(\Omega) \leq H_p(\phi(\Omega))(1 +c  \|(\nabla\phi)^{-1} - I\|_{L^{\infty}(\Omega)})$. Inequality (\ref{bbb}) then follows.
$\hfill\Box$

\section{Differentiability of the Hardy constant}

Let $\Omega $ be a bounded domain   in ${\mathbb{R}}^n$ and  $\psi $  a Lipschitz continuous map  from ${\overline{\Omega}}$ to ${\mathbb{R}}^n$. In the sequel by $\phi_t$ we shall denote the map from $\overline{\Omega} $ to ${\mathbb{R}}^n$
defined by
$$
\phi_t=I +t\psi,
$$
where $t\in {\mathbb{R}}$ and $I$ is the identity map. Clearly, there exists $T>0$ such that for any
$t\in ]-T,T[$ the map $\phi_t$ is a bi-Lipschitz homeomorphism from $\overline{\Omega} $ onto ${\phi_t(\overline{\Omega} )}$,  and  $\phi_t(\partial \Omega)=\partial \phi_t(\Omega )$.

Given a homeomorphism from $\Omega $ onto $\phi (\Omega)$, we set
\begin{equation}\label{viomega}
V_{\Omega }[\phi , \psi](y)=d^{-1}_{ \phi (\Omega )} y
\left( \nabla d_{ \phi (\Omega )}  \right) (y)\cdot\left( \psi\circ \phi^{(-1)}(y)-\psi \circ \phi^{(-1)}(\tau_{ \phi (\Omega )}y) \right),
\end{equation}
for all $y\in \phi (\Omega )$ such that $d_{ \phi (\Omega )}$ is differentiable at $y$. Here and in the sequel, by $\tau_Ax$ we denote the
 nearest point of 
$\partial A$ to $x$, which is unique for almost all $x$.

\begin{lem}\label{der} Let $\Omega $ be a bounded domain in ${\mathbb{R}}^n$, $p\in [1, \infty [$ and $\psi$  a Lipschitz continuous map from $\overline{\Omega} $ to ${\mathbb{R}}^n$. 
Let $T>0$ be such that  $\phi_t=I+t\psi$ is  a  bi-Lipschitz homeomorphism from $\overline{\Omega} $ onto $\phi_t (\overline{\Omega} )$ for all $t\in ]-T,T[$.

 Let $t_0\in ]-T,T[$ be fixed. The following statements hold:
\begin{itemize}
\item[ {\rm (i)} ] There exist  $c,s_0>0$ such that
\begin{equation}\label{der000}
|d^p_{ \phi_{t_0+s} (\Omega )}\phi_{t_0+s}(x)- d^p_{ \phi_{t_0} (\Omega )}\phi_{t_0}(x)  |\le
cd^p_{ \phi_{t_0} (\Omega )}\phi_{t_0}(x)| s|
\end{equation}
for all $x\in \Omega $ and $s\in ]-s_0,s_0[$.
\item[{\rm (ii)} ]
If $x\in\Omega$ and  $d_{ \phi_{t_0} (\Omega )} $ is differentiable at $\phi_{t_0}(x)$ then the map $t\mapsto d^p_{ \phi_{t} (\Omega )}\phi_t(x) $ is differentiable at $t_0$ and
\begin{equation}\label{der0}
\frac{d}{dt}_{|_{t=t_0}} d^p_{ \phi_{t} (\Omega )}\phi_t(x) = pd^{p}_{ \phi_{t_0} (\Omega )} \phi_{t_0}(x)  V_{\Omega }[\phi_{t_0},\psi](\phi_{t_0}(x)).
\end{equation}
\end{itemize}
\end{lem}
{\em Proof.}  It suffices to give a detailed proof only  for the case $p=2$, since the proofs of (\ref{der000}) and (\ref{der0}) for $p\ne 2$ can  be immediately deduced 
from the case $p=2$ combined with inequality (\ref{der3.2}) below. 

Let $x\in \Omega $ and $b\in \partial \Omega $ be such that
$d_{\phi_{t_0}(\Omega)}\phi_{t_0}(x)=|\phi_{t_0}(x)-\phi_{t_0}(b)|$.  For  any $s\in ]-T-t_0, T-t_0[$ we have

\begin{eqnarray}\label{der1} \lefteqn{
d^2_{ \phi_{t_0+s} (\Omega )}\phi_{t_0+s}(x) =\min_{a\in \partial \Omega }| \phi_{t_0+s }(x)-\phi_{t_0+s}(a)  |^2}\nonumber \\
& & \le |\phi_{t_0+s}(x) -\phi_{t_0+s} ( b)  |^2 \nonumber \\
& & = |  \phi_{t_0}(x) - \phi_{t_0}(b)   +s (\psi (x)  -\psi ( b    )  ) |^2\nonumber\\
& & = d^2_{\phi_{t_0}(\Omega )}  \phi_{t_0}(x)+ s^2 |\psi (x)  -\psi (  b  )|^2 \nonumber \\
& & \quad +2s( \phi_{t_0}(x) - \phi_{t_0}(b) )  \cdot ( \psi (x)  -\psi (  b)  ).
\end{eqnarray}

Note that there exists $c>0$ such that
\begin{eqnarray}\label{der1.5}\lefteqn{| \psi (x)  -\psi (  b )|
=|\psi ( \phi_{t_0}^{(-1)} (\phi_{t_0}(x)))  -\psi (  \phi_{t_0}^{(-1)}(  \phi_{t_0}(b) ))|} \nonumber \\
& &\qquad\qquad\qquad\qquad \le c |\phi_{t_0}(x)  - \phi_{t_0}(b)  |=
cd_{ \phi_{t_0} (\Omega )}\phi_{t_0}(x)
\end{eqnarray}
for all $x\in \Omega $. Let $b_{s}\in \partial \Omega$ be such that
\begin{equation}\label{der2}
d_{ \phi_{t_0+s} (\Omega )}\phi_{t_0+s}(x)=|\phi_{t_0+s}(x)-\phi_{t_0+s}(b_s)|.
\end{equation}
Then we have
\begin{eqnarray}\label{der3}
d^2_{ \phi_{t_0+s} (\Omega )}\phi_{t_0+s}(x)&=& |\phi_{t_0}(x)-\phi_{t_0}(b_{s})+ s(\psi(x)-\psi (b_s) )  |^2\nonumber \\
& \geq& d^2_{ \phi_{t_0} (\Omega )}\phi_{t_0}(x)+  2s(  \phi_{t_0}(x)-\phi_{t_0}(b_{s}) )\cdot (\psi(x)-\psi (b_s) )\nonumber \\
& & +s^2|\psi(x)-\psi (b_s) |^2
\end{eqnarray}

Note that there exist $s_0,c>0$ such that
\begin{eqnarray}\label{der3.1} | \psi (x)  -\psi (  b_s )|
&=&|\psi ( \phi_{t_0+s}^{(-1)} (\phi_{t_0+s}(x)))  -\psi (  \phi_{t_0+s}^{(-1)}(  \phi_{t_0+s}(b_s) ))| \nonumber \\
&  \le & c |\phi_{t_0+s}(x)  - \phi_{t_0+s}(b_s)  |  \nonumber \\
&=& cd_{ \phi_{t_0+s} (\Omega )}\phi_{t_0+s}(x)
\end{eqnarray}
for all $x\in \Omega $ and $s\in ]-s_0,s_0[$. Moreover, possibly replacing $s_0$ by a smaller value, there exists $c>0$ such that
\begin{equation}\label{der3.2}
c^{-1}d_{ \phi_{t_0} (\Omega )}\phi_{t_0}(x)\le d_{ \phi_{t_0+s} (\Omega )}\phi_{t_0+s}(x)\le cd_{ \phi_{t_0} (\Omega )}\phi_{t_0}(x),
\end{equation}
for all $x\in\Omega$ and $s\in ]-s_0,s_0[$.
Ineed, the second inequality in (\ref{der3.2}) easily follows by applying (\ref{b6}) with $\Omega$ replaced by $\phi_{t_0}(\Omega)$ and $\phi$ replaced by $\phi_{t_0+s}\circ\phi_{t_0}^{(-1)}$; the first one follows similarly from (\ref{b6}).

From inequalities (\ref{der1})-(\ref{der3.2}) we easily deduce the validity of (\ref{der000}) for $p=2$.

We now assume that $d_{ \phi_{t_0} (\Omega )} $ is differentiable at $\phi_{t_0}(x)$  (hence $\tau_{ \phi_{t_0} (\Omega )} \phi_{t_0}(x)$  is uniquely defined) and prove statement (ii).
By (\ref{der000})  it follows that
\begin{equation}\label{der35}
\lim_{s\to 0}d_{ \phi_{t_0+s} (\Omega )}\phi_{t_0+s}(x)=d_{ \phi_{t_0} (\Omega )}\phi_{t_0}(x).
\end{equation}
We claim that
\begin{equation}\label{der4}
\lim_{s\to 0 } b_s=\phi_{t_0}^{(-1)}(  \tau_{ \phi_{t_0} (\Omega )}\phi_{t_0} (x)  ).
\end{equation}
In order to prove (\ref{der4}) it suffices to prove that
\begin{equation}\label{der5}
\lim_{s\to 0}\phi_{t_0+s}(b_{s})= \lim_{s\to 0 }\phi_{t_0+s}(\phi_{t_0}^{(-1)}(  \tau_{\phi_{t_0} (\Omega )}\phi_{t_0} (x)  )   )
\end{equation}
i.e.,
\begin{equation}\label{der6}
\lim_{s\to 0}\phi_{t_0+s}(b_{s})= \tau_{ \phi_{t_0} (\Omega )}\phi_{t_0} (x)  .
\end{equation}
Assume by contradiction  that (\ref{der6}) doesn't hold. Then there exists $a\in \partial\phi_{t_0}(\Omega )$ such that, possibly passing to a subsequence,
\begin{equation}\label{der7}
\lim_{s\to 0} \phi_{t_0+s}(b_{s})=a \ \ \ {\rm and }\ \  \ | \phi_{t_0}(x) -a|> d_{\phi_{t_0}(\Omega )}\phi_{t_0}(x)+\delta ,
\end{equation}
where $\delta >0$. In particular
\begin{equation}\label{der8}
\lim_{s\to 0}| \phi_{t_0+s}(b_s)-\phi_{t_0}(x)  |=|a-\phi_{t_0}(x) |>  d_{\phi_{t_0}(\Omega )}\phi_{t_0}(x)+\delta .
\end{equation}
We also have
\begin{eqnarray}\label{der9} \lefteqn{|  \phi_{t_0+s}(b_s)-\phi_{t_0}(x)  |^2= | \phi_{t_0+s}(b_s)-\phi_{t_0+s}(x)+s\psi (x)    |^2}\nonumber \\
& &
=d^2_{\phi_{t_0+s}(\Omega)}\phi_{t_0+s}(x)+2s( \phi_{t_0+s}(b_s)-\phi_{t_0+s}(x)  )\cdot \psi (x) +s^2|\psi (x)|^2 .
\end{eqnarray}
By (\ref{der35}) and (\ref{der9}) we deduce that
$$
\lim_{s\to 0}| \phi_{t_0+s}(b_s)-\phi_{t_0}(x) |= d_{\phi_{t_0}(\Omega )}\phi_{t_0}(x)
$$
which contradicts (\ref{der8}). Thus (\ref{der4}) holds.

By (\ref{der1}), (\ref{der3}) and (\ref{der4}), by observing that $b=\phi_{t_0}^{(-1)}(\tau_{ \phi_{t_0} (\Omega )}\phi_{t_0} (x)  ) $ and noting that 
$$
(\nabla d_{\phi_{t_0}(\Omega )})(\phi_{t_0}(x) )=\frac{\phi_{t_0}(x) - \tau_{ \phi_{t_0} (\Omega )}\phi_{t_0} (x)}{d_{\phi_{t_0}(\Omega)}\phi_{t_0}(x)  },
$$
we immediately deduce the validity of (\ref{der0}) for $p=2$. 
\hfill $\Box$

\begin{lem}\label{den}
Let $\Omega $ be a bounded domain in ${\mathbb{R}}^n$, $p\in [1, \infty [$ and $\psi$  a Lipschitz continuous map from $\overline{\Omega} $ to ${\mathbb{R}}^n$. 
Let $T>0$ be such that  $\phi_t=I+t\psi$ is  a  bi-Lipschitz homeomorphism from $\overline{\Omega} $ onto $\phi_t (\overline{\Omega} )$ for all $t\in ]-T,T[$.
 Let $u\in W^{1,p}_0(\Omega )$, $\rho\in L^{\infty}(\Omega)$. Then the function  $G :]-T,T[ \to {\mathbb{R}}$
defined by
\begin{equation}\label{den0}
G(t)=\int_{\Omega }\frac{|u|^p\rho }{d^p_{\phi_t(\Omega )}\phi_t(x)}dx,
\end{equation}
for all $t\in ]-T,T[$,  is differentiable and
\begin{equation}\label{den1}
G'(t)=-p\int_{\Omega } \frac{|u|^p\rho  V_{\Omega }[\phi_t,\psi ](\phi_t(x))}{d^p_{ \phi_t(\Omega )}\phi_t(x)}dx,
\end{equation}
for all $t\in ]-T,T[$.
\end{lem}
{\em Proof.} We divide the proof into two steps.

Step 1. We prove that $G$ is differentiable at $t=0$ and that formula (\ref{den1}) holds for $t=0$.  Obviously, we have
\begin{equation}\label{den6}
\frac{G(t)-G(0)}{t}=-\int_{\Omega }\frac{|u|^p\rho  (d^p_{\phi_t(\Omega)}\phi_t(x)-d^p_{\Omega}x) }{td^p_{\phi_t(\Omega)}\phi_t(x)d^p_{\Omega}x}dx.
\end{equation}
By (\ref{der000}) we have that there exists $c>0$ and $T_0\in ]0,T[$  such that
$$
\frac{|u|^p\rho  |d^p_{\phi_t(\Omega)}\phi_t(x)-d^p_{\Omega}x| }{|t|d^p_{\phi_t(\Omega)}\phi_t(x)d^p_{\Omega}x}
\le c \frac{|u|^p\rho}{d^p_{\Omega}x },$$
for all $t\in ]-T_0,T_0[$ and $x\in\Omega$. Since the function $|u|^p\rho /d^p_{ \Omega }$ belongs to $L^1(\Omega )$ and does not depend on $t\in ]-T_0,T_0[$, we
can apply the  Dominated Convergence Theorem and  pass to the limit under the integral sign in (\ref{den6}) as $t\to 0$. By applying formula (\ref{der0}) with $t_0=0$ we immediately get that $G$ is differentiable at $t=0$ and that formula (\ref{den1}) holds  for $t=0$.

Step 2. Let $t_0\in ]-T,T[$ be fixed. By changing variables in integrals, we have
\begin{equation}\label{den70}
G(t)= \int_{\phi_{t_0}(\Omega )}\frac{(|u|^p\rho)\circ\phi_{t_0}^{(-1)}|{\rm det }\nabla \phi_{t_0}^{(-1)} |   }{d^p_{ \phi_t\circ\phi_{t_0}^{(-1)} (\phi_{t_0}(\Omega ))}(\phi_t\circ\phi_{t_0}^{(-1)}(y))  }dy.
\end{equation}

We set $t=t_0+s$ and we note that
\begin{eqnarray}\label{den8}
\phi_t\circ \phi_{t_0}^{(-1)}(y) = \phi_{t_0}^{(-1)}(y)+(t_0+s)\psi(\phi_{t_0}^{(-1)}(y))=y+s\psi (\phi_{t_0}^{(-1)}(y)),
\end{eqnarray}
for all $y\in \phi_{t_0}(\Omega)$. By setting $\tilde \Omega =\phi_{t_0}(\Omega )$,
$\tilde u =u\circ \phi_{t_0}^{(-1)}$, $\tilde \rho =\rho\circ \phi_{t_0}^{(-1)} |{\rm det }\nabla \phi_{t_0}^{(-1)} |$,  $\tilde \psi =\psi\circ \phi_{t_0}^{(-1)}$,
$\tilde \phi_s =I+s\tilde \psi$,
we have
$$
G(t_0+s)=\int_{\tilde \Omega }\frac{| \tilde u|^p\tilde \rho }{d^p_{\tilde\phi_s (\tilde \Omega ) }\tilde \phi_s(y)}dy.
$$
By Step 1, it follows that $G$ is differentiable at $t_0$ and
\begin{equation}\label{den9}
G'(t_0)=-p\int_{\tilde \Omega } \frac{| \tilde u|^p\tilde\rho V_{\tilde\Omega }[\tilde\phi_0,\tilde \psi ](\tilde \phi_0 (y))}{d^p_{ \tilde \phi_0(\tilde\Omega )}\tilde\phi_0(y)}dy.
\end{equation}
It is easily seen that
\begin{equation}\label{den7}
V_{\tilde\Omega }[\tilde\phi_0,\tilde \psi ](y  )=V_{\Omega }[\phi_{t_0},\psi ](y)
\end{equation}
By changing variables in the right-hand side of (\ref{den9}) and using (\ref{den7}) we get formula (\ref{den1}) for $t=t_0$. \hfill $\Box$

In order to prove  that the Hardy constant $H_p(\phi_t (\Omega))$ is differentiable with respect to $t$ in the case  $H_p(\phi_t (\Omega))<((p-1)/p)^p$, we need to prove a result concerning the continuous dependence on $t$ of the corresponding minimizers. To do so, we need the following  theorem which provides  estimates for the minimizers and their gradients.  As we have already mentioned before, estimate (\ref{marshaf1}) is proved in \cite{MMP, mash}. Here we indicate the dependence of  the constant $C$ in (\ref{marshaf1}) on the data, and we prove estimate (\ref{marshaf1grad}) which is also  of independent interest.

For $\delta>0$ we define
\[
\Omega_{\delta}:=\{x\in\Omega : d_{\Omega}(x)<\delta\} .
\]
By ${\rm Inr}(\Omega )$ we denote the inradius of $\Omega$,  namely ${\rm Inr}(\Omega )=\sup_{\Omega}d_{\Omega}$.
Moreover,   by $\| v\|_{L^p(\Omega , \rho)}$ we denote the weighted norm $(\int_{\Omega}|v|^p\rho dx)^{1/p}$ of a function $v$ defined on $\Omega$.

\begin{thm}\label{marshaf} Let $p\in ]1,\infty [$ and  let $\Omega$ be a bounded domain in ${\mathbb{R}}^n$ of class $C^2$ such that $H_p(\Omega )<((p-1)/p)^p$.  Let $\delta>0$ be such that $d_{\Omega}$ is of class $C^2$  on $\Omega_{2\delta}$. Let
$C_{\delta}>0$ be such that 
\begin{equation}\label{marshaf0}
|\Delta d_{\Omega}|\le C_{\delta},\ \ \ {\rm on}\ \Omega_{\delta}.
\end{equation}
Then there exists $C>0$ which depends only on $n,p, \delta,C_{\delta}$ and ${\rm Inr}(\Omega) $ such that if $v$ is a positive minimizer for
the Hardy constant $H_p(\Omega)$ normalized by $\|v\|_{L^p(\Omega ,d_{\Omega}^{-p})}=1 \, $, then
\begin{equation}\label{marshaf1}
v\le C d^{\alpha}_{\Omega},\ \ \ {\rm on}\ \Omega, 
\end{equation}
and 
\begin{equation}\label{marshaf1grad}
|\nabla v|\le C d^{\alpha-1}_{\Omega},\ \ \ {\rm on}\ \Omega, 
\end{equation}
where $\alpha$  is the largest solution to the equation $(p-1)\alpha^{p-1}(1-\alpha )=H_p(\Omega)$. 
\end{thm}

{\em Proof.}   The existence of a constant $C>0$ such that inequality (\ref{marshaf1}) holds, is proved in \cite[Lemma~9]{MMP} and \cite[Lemma~5.2]{mash}. Moreover, a detailed analysis of the proof in \cite{mash} combined with the local estimates in Serrin~\cite[Theorems 1,~2]{se} allow to deduce that the constant $C$ in (\ref{marshaf1}) depends only on  $n,p, \delta,C_{\delta}$ and ${\rm Inr}(\Omega) $.   

 We now prove (\ref{marshaf1grad}).
Let $x_0\in \Omega $ and $R>0$ be such that $B(x_0,2R)\subset \Omega$.  Note that $v$ is a solution to equation (\ref{pde}). Thus by the general gradient
estimate in \cite[Theorem 1.1]{dumi} there exists $c>0$  such that 
\begin{equation} \label{ming1}
|\nabla v(x_0)|\le c\avint_{B(x_0,R)}|\nabla v|dx+c\int_0^R\left(\frac{|\mu |(B(x_0,\rho))}{\rho^{n-1}}\right)^{\frac{1}{p-1}}\frac{d\rho}{\rho},
\end{equation}
where $\mu$ is the measure with density  $ H_p(\Omega )v ^{p-1}/d^p_{\Omega }$, hence
$|\mu |(B(x_0,R))$ equals  $ H_p(\Omega )\int_{B(x_0,\rho)}v^{p-1}/d_{\Omega}^pdx$. (Note that in our case it suffices to have $\mu \in L^1_{loc}(\Omega)$ in order to apply  \cite[Theorem 1.1]{dumi}.)
We set $R=d_{\Omega}(x_0)/3$ and we observe that for  any $x\in B(x_0,R)$ we have that 
\begin{equation}\label{ming3}
R\le d_{\Omega}(x)\le 4R.
\end{equation}
We now estimate the first summand in the right-hand sinde of (\ref{ming1}).  
Note that by \cite[Theorems 1,2]{se} we have that 
\begin{equation}\label{serringrad}
\| \nabla u\|_{L^p(B(x_0,R))}\le \frac{c}{R}\| u\|_{L^p(B(x_0,2R))}.
\end{equation}
Thus, by H\"{o}lder's inequality, (\ref{marshaf1}), (\ref{ming3}) and (\ref{serringrad}) we get
\begin{eqnarray}\lefteqn{
\avint_{B(x_0,R)}|\nabla v|dx\le c\frac{|B(x_0,R)|^{-\frac{1}{p}}}{R}\| u\|_{L^p(B(x_0,2R))}\le cR^{-\frac{n}{p}-1}\| d^{\alpha}_{\Omega}\|_{L^p(B(x_0,2R))}}\nonumber \\ & &\qquad \qquad\quad
\le cR^{-\frac{n}{p}-1}\left(\int_{B(x_0,2R)}R^{\alpha p}dx\right)^{\frac{1}{p}}\le cR^{\alpha-1}\le c d^{\alpha-1}_{\Omega}(x_0).
\label{mingpre5}
\end{eqnarray}
We now estimate the second summand in the right-hand sinde of (\ref{ming1}).  
By (\ref{marshaf1}) and (\ref{ming3}) it follows that
\begin{equation}
\label{ming4}
|\mu |(B(x_0,\rho))|\le c\int_{B(x_0,\rho)}d_{\Omega}^{\alpha (p-1)-p}dx\le c\int_{B(x_0,\rho)}R^{\alpha (p-1)-p}dx\le cR^{\alpha (p-1)-p}\rho^n,
\end{equation}
for all $\rho \in ]0,R[$.
By (\ref{ming4}) we immediately get
\begin{equation}
\label{ming5}
\int_0^R\left(\frac{|\mu |(B(x_0,\rho))}{\rho^{n-1}}\right)^{\frac{1}{p-1}}\frac{d\rho}{\rho}\le c R^{\alpha -1}\le c d^{\alpha-1}_{\Omega}(x_0).
\end{equation}
Inequalities (\ref{ming1}), (\ref{mingpre5}) and (\ref{ming5}) imply the validity of (\ref{marshaf1grad}). \hfill $\Box$

By $C^2(\overline{\Omega};{\mathbb{R}}^n)$ we denote the space of functions $\phi$ from $\overline{\Omega} $ to ${\mathbb{R}}^n$ of class $C^2$ in $\Omega$  such that all derivatives $D^{\alpha}\phi $ with $0\le |\alpha| \le 2 $ have continuous extensions on $\overline{\Omega}$. Moreover, we endow it with its standard
norm $\| \phi \|_{C^2(\overline{\Omega} ;{\mathbb{R}}^n)}=\max_{0\le |\alpha |\le 2}\max_{x\in \overline{\Omega}} |D^{\alpha}\phi (x)|$.

\begin{corol} \label{marcust0}
 Let $\Omega$ be a bounded domain in ${\mathbb{R}}^n$ of class $C^2$ and $M>0$. Let  $\psi\in C^2(\overline{\Omega} ;{\mathbb{R}}^n)$ with $\| \psi \|_{C^2(\overline{\Omega} ;{\mathbb{R}}^n)}\le M$. Let $T>0$ be such that  $\phi_t=I+t\psi$ is  a diffeomorphism from $\overline{\Omega} $ onto $\phi_t (\overline{\Omega} )$ for all $t\in ]-T,T[$. 

Then there exist $T_0\in ]0,T[, R, C_R>0$ depending only on $\Omega $ and $M$ such that  $d_{\phi_t (\Omega)}$ is of class $C^2$ on $(\phi_t(\Omega ))_{R}$ and 
\begin{equation}\label{gilba}|\Delta d_{\phi_t (\Omega)} |\le C_R,\ \ \ {\rm on}\ \ (\phi_t(\Omega))_{R},\end{equation}
for all $t\in ]-T_0,T_0[$.
Furthermore, if $p\in ]1,\infty [$ and $H_p(\Omega)<((p-1)/p)^p$ then there exist $T_1\in ]0,T_0[$  and   $C>0$  depending only on $\Omega, p $ and $M$, such that   
$H_p(\phi_t(\Omega ))<((p-1)/p)^p$  for all $t\in ]-T_1,T_1[$  and such that if $v_t$ is a positive minimizer for $H_p(\phi_t(\Omega ))$ normalized by $\|v_t\|_{L^p(\phi_t(\Omega ), d_{\phi_t(\Omega )}^{-p})} =1$  then
\begin{equation}\label{marcust1}
v_t\le C d_{\phi_t(\Omega )}^{\alpha_t }\ \ {\rm and}\ \ |\nabla v_t|\le C d_{\phi_t(\Omega )}^{\alpha_t -1},\ \ {\rm on}\ \ \phi_t(\Omega ),
\end{equation}
for all $t\in ]-T_1,T_1[$,
where $\alpha_t $ is the largest solution to the equation $(p-1)\alpha_t^{p-1}(1-\alpha_t )=H_p(\phi_t(\Omega) )$.
\end{corol}

{\em Proof.} The first part of the statement follows by standard arguments. In particular,  we refer to Gilbarg and Trudinger~\cite[Appendix~14.6]{gitr}
for the proof of (\ref{gilba}). Assume now that $H_p(\Omega)<((p-1)/p)^p$. By Proposition \ref{thm:bbb} it follows that $H_p(\phi_t(\Omega ))<((p-1)/p)^p$  for all $t$ sufficiently small. The proof of (\ref{marcust1}) immediately follows by Theorem \ref{marshaf} and (\ref{gilba}).\hfill $\Box$

We are now ready to prove the following 

\begin{thm}[Stability of minimizers]  
\label{stabmin} Let $p\in ]1,\infty [$ and  $\Omega$ be a bounded domain in ${\mathbb{R}}^n$ of class $C^2$ such that  $H_p(\Omega)<((p-1)/p)^p$.
Let $\psi\in C^2(\overline{\Omega};{\mathbb{R}}^n)$ and let $T>0$ be such that  $\phi_t=I+t\psi$ is a  diffeomorphism from $\overline{\Omega} $ onto $\phi_t (\overline{\Omega} )$ for all $t\in ]-T,T[$. 

Let $T_0\in ]0,T[$ be such that $H_p(\phi_t (\Omega ))<((p-1)/p)^p$ for all $t\in ]-T_0,T_0[$ and let $v_t$ be a positive minimizer for $H_p(\phi_t(\Omega ))$ normalized by $\|v_t\|_{L^p(\phi_t(\Omega ), d^{-p}_{\phi_t(\Omega )})} =1$.  
Then the following statements hold:
\begin{itemize}
\item[ { (i) } ] Let $u_t=v_t\circ \phi_t$ for all $t\in ]-T_0,T_0[$. Then 
\[ \lim_{t\to 0}\| u_t-u_0\|_{W^{1,p}_0(\Omega)}=0.\]
\item[{ (ii) } ] Assuming that every function $v_t$  is extended by zero outside $\phi_t(\Omega)$,  we have 
\[ \lim_{t\to 0}\| v_t-v_0\|_{W^{1,p}_0(\Omega\cup \phi_t (\Omega))}=0 .\]
\end{itemize}
\end{thm}
{\em Proof.}  First, we prove statement (i). By the normalization of $v_t$ we have that $\| \nabla v_t\|^p_{L^p(\phi_t(\Omega ))}=H_p(\phi_t(\Omega))$. It follows from Proposition
\ref{thm:bbb} that $\| v_t\|_{W^{1,p}_0(\phi_t(\Omega))}$ and
$\| u_t\|_{W^{1,p}_0(\Omega)}$ are uniformly bounded for $t$ sufficiently small. Thus, possibly passing to a subsequence, there exists $\tilde u_0\in W^{1,p}_0(\Omega)$ such that $u_t\to \tilde u_0$  weakly in $W^{1,p}_0(\Omega)$ and strongly in $L^p(\Omega)$ as $t\to 0$. 

We claim that 
\begin{equation}\label{stabmin05}
\| \tilde u_0\|_{L^p(\Omega , d_{\Omega}^{-p})}=1.
\end{equation}
 In order to prove this, we recall first that 
\begin{equation}\label{stabmin1}
\int_{\Omega} \frac{u_t^p(x) |{\rm det}\nabla \phi_t(x)|}{d^p_{\phi_t(\Omega )}\phi_t(x)}dx=1
\end{equation}
for all $t\in ]-T,T[$. By the continuous dependence of $H_{\phi_t(\Omega)}$ on $t$ proved in Proposition \ref{thm:bbb},  and by Corollary~\ref{marcust0}, there exist $C>0$ and $\alpha >(p-1)/p$ independent of $t$ such that 
\begin{equation}\label{stabmin2}
u_t(x)\le C d^{\alpha}_{\phi_t(\Omega )}\phi_t(x)
\end{equation} 
for all $x\in \Omega$ and $t$ sufficiently small. By (\ref{der000}) we have that $d^p_{\phi_t(\Omega )}\phi_t(x)\to d^p_{\Omega }(x)$ as $t\to 0$. Thus by  (\ref{der000}),  (\ref{stabmin2}) and the Dominated Convergence Theorem we can pass  to the limit inside the integral sign in (\ref{stabmin1}) and prove the claim above.

By the elementary inequality 
$
|a|^p\geq |b|^p+p|b|^{p-2}b\cdot (a-b)
$
valid for  all vectors $a, b$ in ${\mathbb{R}}^N$,  we have
\begin{eqnarray}
& &H_p(\phi_t(\Omega )) =\int_{\Omega }|\nabla u_t(\nabla \phi_t)^{-1}|^{p}|{\rm det}\nabla \phi_t(x)| dx\geq \int_{\Omega}|\nabla \tilde u_0|^p|{\rm det}\nabla \phi_t(x)| dx\nonumber  \\
& & +p\int_{\Omega }|\nabla \tilde u_0|^{p-2}\nabla \tilde u_0\cdot ( \nabla u_t(\nabla \phi_t)^{-1}-\nabla \tilde u_0  )|{\rm det}\nabla \phi_t(x)|dx.
\end{eqnarray}
By passing to the limit in the previous inequality and using Proposition~\ref{thm:bbb}, we get
$$
H_p(\Omega )\geq \int_{\Omega }|\nabla \tilde u_0|^{p}.
$$
This, combined with the normalization of $\tilde u_0$ and its positivity, allows to conclude that  $\tilde u_0=u_0$.
Since
\begin{equation}\label{stabmin5}
H_p(\phi_t(\Omega))-\int_{\Omega}|\nabla u_t|^pdx=\int_{\Omega }|\nabla u_t(\nabla\phi_t)^{-1}|^p|{\rm det}\nabla \phi_t(x)| -|\nabla u_t|^pdx ,\\
\end{equation}
 we have that 
\begin{equation}\label{stabmin6}\lim_{t\to 0} \Big(
H_p(\phi_t(\Omega))-\int_{\Omega}|\nabla u_t|^pdx\Big)=0.
\end{equation}
Thus from (\ref{stabmin6}) and  Proposition \ref{thm:bbb} we deduce that 
\begin{equation}\label{stabmin7}\lim_{t\to 0} 
\int_{\Omega}|\nabla u_t|^pdx=H_p(\Omega )=\int_{\Omega }|\nabla u_0|^pdx.
\end{equation}
By (\ref{stabmin7}) and the weak convergence of $u_t$ to $u_0$ in $W^{1,p}_0(\Omega )$  we deduce the validity of (i).

We now prove statement (ii). We have
\begin{eqnarray}\lefteqn{
\| \nabla v_t -\nabla v_0\|^p_{L^p(\Omega\cup\phi_t(\Omega))} } \nonumber  \\ & & 
= \|\nabla v_0\|^p_{L^p( \Omega \setminus \phi_t(\Omega) )} 
+ \| \nabla v_t \|^p_{L^p(\phi_t(\Omega)\setminus \Omega)} +
\| \nabla v_t -\nabla v_0\|^p_{L^p(\Omega\cap  \phi_t(\Omega))}.
\label{tz}
\end{eqnarray}
The first summand in (\ref{tz}) above clearly tends to zero as $t\to 0$. The fact that the second summand in (\ref{tz}) also tends to zero will follow immediately from
the following

{\bf Claim.} For all $\epsilon>0$ there exist $ \tau , \delta>0$ such that
$$\int_A|\nabla v_t|^pdy <\epsilon , $$
for all $t\in ]-\tau, \tau[$ and  all measurable subsets $A$ of $ \phi_t(\Omega)$ with $|A|<\delta$.

To prove the Claim we use   (\ref{der000}) and Corollary \ref{marcust0} to conclude that there exists
$\alpha_* \in ](p-1)/p,\alpha[$ such that
$$
\| \nabla v_t \|^p_{L^p(A)} \leq  c\int_{A} 
d_{\phi_t(\Omega)}^{p(\alpha_t-1)} dy 
\leq c \int_{\phi_t^{(-1)}(A)} d_{\Omega}^{p(\alpha_t-1)} dx 
\leq c \int_{\phi_t^{(-1)}(A)} d_{\Omega}^{p(\alpha_*-1)} dx ,
$$
for all small enough $t$; the last integral clearly tends to zero uniformly with respect to $t$ as $|A|\to 0$.

To estimate the third summand in (\ref{tz}) we take a set $U\subset\subset\Omega$ and write
\begin{eqnarray}\lefteqn{
 \| \nabla v_t -\nabla v_0\|^p_{L^p(\Omega\cap  \phi_t(\Omega))} \leq \| \nabla v_t -\nabla v_0\|^p_{L^p(U\cap  \phi_t(U))}}\nonumber   \\ 
& & \qquad  +  \| \nabla v_t -\nabla v_0\|^p_{L^p( (\Omega\setminus U) \cap  \phi_t(\Omega))}
  + \| \nabla v_t -\nabla v_0\|^p_{L^p(\Omega\cap  \phi_t(\Omega\setminus U))}.\label{tz005}
\end{eqnarray}
By the last Claim, the last two norms in (\ref{tz005}) can be made aritrarily small provided $U$ is a large enough subset of $\Omega$.  Hence the proof will be complete if we show that for a fixed $U\subset\subset\Omega$ there holds
\begin{equation}
\| \nabla v_t -\nabla v_0\|^p_{L^p(U\cap \phi_t(U))} \longrightarrow 0, \; \qquad \mbox{ as } t\to 0.
\label{tz1}
\end{equation}
Indeed we have
\begin{eqnarray}
&&\hspace{-0.4cm} \| \nabla v_t -\nabla v_0\|_{L^p(U\cap  \phi_t(U)   )} \nonumber \\
&=& \| [(\nabla u_t)(\nabla\phi_t)^{-1}] \circ \phi_t^{(-1)}   -\nabla v_0\|_{L^p(U\cap  \phi_t(U) )}\nonumber \\
&\leq& c \| [(\nabla u_t)(\nabla\phi_t)^{-1}]   -(\nabla v_0)\circ \phi_t\|_{L^p(\phi_t^{(-1)}(U)\cap  U )}\nonumber \\
&\leq& c \| [(\nabla u_t)(\nabla\phi_t)^{-1}]  -\nabla u_t \|_{L^p(\phi_t^{(-1)}(U)\cap  U )}
+  c \| \nabla u_t -\nabla u_0 \|_{L^p(\phi_t^{(-1)}(U)\cap  U )} \nonumber \\
&& +   c \| \nabla v_0 -(\nabla v_0)\circ\phi_t \|_{L^p(\phi_t^{(-1)}(U)\cap  U )} .
\label{tz2}\end{eqnarray}
The first of the last three terms in (\ref{tz2}) clearly tends to zero as $t\to 0$. The same is true for the second term by statement (i).  The fact that the third term tends to zero is an immediate consequence of the  local H\"{o}lder continuity of $\nabla v_0$,  which follows by the general results in \cite{DiBe}. This completes the proof.\hfill $\Box$

Finally, we can prove the following Hadamard-type formula for the $L^p$ Hardy constant.  Note  that  by Corollary~\ref{marcust0}, the assumption $H_p(\Omega )<((p-1)/p)^p$ guarantees the existence of $T_0\in ]0,T[$ such that $H_p(\phi_t(\Omega ))<((p-1)/p)^p$ for all $t\in ]-T_0,T_0[$.

\begin{thm}
\label{had}Let $p\in ]1,\infty [$ and $\Omega$ be a bounded  domain in ${\mathbb{R}}^n$ of class $C^2$ such that $H_p(\Omega )<((p-1)/p)^p$. Let $\psi\in C^2(\overline{\Omega};{\mathbb{R}}^n)$ and let $T>0$ be such that  $\phi_t=I+t\psi$ is a  diffeomorphism from $\overline{\Omega} $ onto $\phi_t (\overline{\Omega} )$ for all $t\in ]-T,T[$. Let $T_0\in ]0,T[$ be such that $H_p(\phi_t(\Omega ))<((p-1)/p)^p$  for all $t\in ]-T_0,T_0[$. 

Then  $H_p(\phi_t(\Omega ))$ is differentiable with respect to $t$ for all $ t\in ]-T_0,T_0[$. Moreover
\begin{eqnarray}\label{had1}\lefteqn{
\frac{dH_p(\phi_t(\Omega ))}{dt} = \int_{\phi_t(\Omega ) }|\nabla v_t|^p{\rm div }(\psi \circ\phi_t^{(-1)}) -p|\nabla v_t|^{p-2}\nabla v_t\nabla (\psi \circ\phi_t^{(-1)}) (\nabla v_t)^t dy } \nonumber \\
& &\qquad\quad+ H_p(\phi_t(\Omega ) )\int_{\phi_t(\Omega ) } \frac{v_t^p}{d^p_{\phi_t(\Omega ) }}\left(  pV_{\Omega }[\phi_t,\psi ] -  {\rm div}(\psi\circ\phi_t^{(-1)}) \right)dy,
\end{eqnarray}
for all $t\in ]-T_0,T_0[$, where $v_t$ is a minimizer for $H_p(\phi_t(\Omega ))$ normalized by the condition $\| v_t\|_{L^p(\phi_t(\Omega) ,  d^{-p}_{\phi_t(\Omega)})}  =1 $.
\end{thm}
{\em Proof.} First of all we note that it suffices to prove that the map $t\mapsto H_p(\phi_t(\Omega ))$ is differentiable
 at $t=0$ and that formula (\ref{had1}) holds for $t=0$. Indeed, as in the proof of Lemma~\ref{den}, if $t\ne 0$ one can consider $\phi_{t}(\Omega )$ as a reference domain subject to the domain transformations
 $\tilde \phi_s=I+s\psi\circ \phi_{t}^{(-1)}$, $s\in {\mathbb{R}}$, and apply the formula for  $s=0$. 
  
 Let $v_t\in W^{1,p}_0(\phi_t(\Omega ))$ be a positive minimizer for $H_p(\phi_t(\Omega ))$ normalized as in the statement.  
We then have (cf. (\ref{chvar})),
\begin{equation}\label{had1bis}
  H_p(\phi_t(\Omega ))=\min_{\substack{u\in W^{1,p}_0(\Omega )\\ u\ne 0 }} R_t[u],
  \end{equation}
  where $R_t[u]=N_t[u]/D_t[u]$ and 
 \begin{equation}\label{had2}
 N_t[u]=\int_{\Omega }|\nabla u(\nabla\phi_t)^{-1}|^p|{\rm det \nabla \phi_t}|dx,
 \end{equation}
 \begin{equation}\label{had3}
 D_t[u]=\int_{\Omega } \frac{u^p}{d^p_{\phi_t(\Omega )}\phi_t(x) }   |{\rm det \nabla \phi_t}|dx. 
 \end{equation}
 We set $u_t=v_t\circ \phi_t$. Clearly, $u_t$ is a minimizer in (\ref{had1bis}) and $H_p(\phi_t(\Omega ))=R_t[u_t]=N_t[u_t]$.
 By Lemma~\ref{den} and standard calculus it follows that the functions $N_t[u]$ and $D_t[u]$ are differentiable with respect to $t$ on $ ]-T_0, T_0[$. By definition, it follows that 
 \begin{equation}\label{had4}
 R_t[u_t]-R_0[u_t]\le H_p(\phi_t(\Omega ))- H_p(\phi_0(\Omega ))\le R_t[u_0]-R_0[u_0],
 \end{equation}
 hence  by the Mean Value Theorem, it follows that  there exist real numbers $\xi (t)$, $\eta (t) $ with $|\xi (t)|, |\eta (t)|< |t|$, such that
 \begin{equation}\label{had5}
 R_{\xi (t)}'[u_t]t\le H_p(\phi_t(\Omega ))- H_p(\phi_0(\Omega ))\le R'_{\eta (t)}[u_0]t,
 \end{equation}
where by $R'_a$ we denote the partial derivative of $R_t[u]$ with respect to $t$ at the point $t=a$ (the same notation is used below for the derivatives of $N$ and $D$).
By standard calculus, we have that for any $u\in H^1_0(\Omega )$,
\begin{equation}
\label{had6}
\frac{d}{dt}|\nabla u(\nabla\phi_t)^{-1}|^p   = -p |\nabla u(\nabla \phi_t )^{-1}|^{p-2} \nabla u 
  (\nabla\phi_t)^{-1} \nabla \psi (\nabla \phi_t)^{-1}  (\nabla\phi_t^{-t}) (\nabla u)^t
\end{equation}
and
\begin{equation}\label{had7}
\frac{d}{dt}|{\rm det}\nabla \phi_t|=\left( \frac{{\rm div }(\psi\circ\phi_t^{(-1)})  }{| {\rm det}\nabla \phi_t^{(-1)}| } \right)\circ \phi_t \, .
\end{equation}
Thus, by (\ref{had6}), (\ref{had7}) and Lemma~\ref{den} we have that for any $u\in H^1_0(\Omega )$
\begin{eqnarray}\lefteqn{
N'_t[u]=\int_{\Omega }|\nabla u(\nabla\phi_t)^{-1}|^p  \left( \frac{{\rm div }(\psi\circ\phi_t^{(-1)})  }{| {\rm det}\nabla \phi_t^{(-1)}| } \right)\circ \phi_t \,   dx} \label{had8}  \\
& & \qquad \qquad  - p \int_{\Omega }  |\nabla u(\nabla \phi_t )^{-1}|^{p-2}  \nabla u  (\nabla\phi_t)^{-1} \nabla \psi  (\nabla\phi_t)^{-t} (\nabla u)^t  |{\rm det \nabla \phi_t}|dx  
\nonumber
\end{eqnarray}
and
\begin{eqnarray}
\lefteqn{
D'_t[u]=\int_{\Omega }  \frac{u^p}{d^p_{\phi_t(\Omega )}\phi_t(x) }\left( \frac{{\rm div }(\psi\circ\phi_t^{(-1)})  }{| {\rm det}\nabla \phi_t^{(-1)}| } \right)\circ \phi_t   \, dx}\nonumber    \\
& & \qquad\qquad   -p\int_{\Omega } \frac{u^p    V_{\Omega }[\phi_t,\psi ](\phi_t(x))}{d^p_{\phi_t(\Omega )}\phi_t(x)}|{\rm det \nabla \phi_t}|dx\label{had9}
\end{eqnarray}
By (\ref{had8}) and (\ref{had9}) and the Dominated Convergence Theorem, it follows that 
\begin{equation}\label{had10}
\lim_{t\to 0} R'_{\eta (t)}[u_0]= R'_{0}[u_0].
\end{equation}
Similarly, and using also Theorem~\ref{stabmin}, we have
\begin{equation}\label{had11}
\lim_{t\to 0} R'_{\xi (t)}[u_t]= R'_{0}[u_0].
\end{equation}
From (\ref{had5}), (\ref{had10}) and (\ref{had11}) it immediately follows that $H_p(\phi_t (\Omega))$ is differentiable with respect to $t$ at $t=0$
and, taking into account that $u_0$ is normalized, we get
\begin{eqnarray}
\label{had12}
\frac{dH_p(\phi_t(\Omega ))}{dt}_{|t=0}&=&\frac{N'_0[u_0]}{D_0[u_0]} -\frac{N_0[u_0]D'_0[u_0]}{D_0^2[u_0]} \\
&= & N'_0[u_0]-H_p(\Omega )D'_0[u_0] \nonumber \\
& =& \int_{\Omega }|\nabla u_0|^p{\rm div }\psi \, dx-p\int_{\Omega }|\nabla u_0|^{p-2}\nabla u_0\nabla\psi (\nabla u_0)^t dx \nonumber \\
& &  -H_p(\Omega )\int_{\Omega }\frac{u_0^p{\rm div}\psi}{d^p_{\Omega }}dx  +pH_p(\Omega )\int_{\Omega }\frac{u_0^p}{d^p_{\Omega }}V_{\Omega }[\phi_0,\psi ]dx  \nonumber
\end{eqnarray}
as required. \hfill $\Box$

Combining Theorems \ref{stabmin} and  \ref{had} we can prove a Fr\'{e}chet differentiability result. Namely, given a bounded open set $\Omega $, we set
\begin{equation}
{\mathcal{A}}_{\Omega}=\Big\{\phi \in C^2(\overline{\Omega} ;{\mathbb{R}}^n):\ \min_{\overline{\Omega}}|{\rm det }\nabla \phi|>0\Big\}
\end{equation}
We then have
\begin{thm}
\label{thm:diff}
Let   $\Omega$ be a bounded domain in ${\mathbb{R}}^n$ of class $C^2$ and $p\in ]1,\infty [$. Then the set
\begin{equation}
{\mathcal{H}}_{\Omega}=\biggl\{\phi \in {\mathcal{A}}_{\Omega}:\ H_p(\phi (\Omega))<\Big(\frac{p-1}{p}\Big)^p\biggr\}
\end{equation}
is open in the space $C^2(\overline{\Omega} ;{\mathbb{R}}^n)$, and the map $H$ from ${\mathcal{H}}_{\Omega}$ to ${\mathbb{R}}$ which takes any $\phi\in {\mathcal{H}}_{\Omega}$ to $ H_p(\phi (\Omega ))$ is Fr\'{e}chet differentiable.   Moreover, the Fr\'{e}chet differential of $H$ at a point $\phi \in {\mathcal{H}}_{\Omega}$ is given by the formula
\begin{eqnarray}
\label{frech}
&& dH_p(\phi )[\psi] = \int_{\phi (\Omega ) }\big\{ |\nabla v|^p{\rm div }(\psi \circ\phi^{(-1)}) -p|\nabla v|^{p-2}\nabla v\nabla (\psi \circ\phi^{(-1)})  (\nabla v)^t \big\} dy  \nonumber \\
&&\qquad\qquad \quad + H_p(\phi(\Omega ) )\int_{\phi (\Omega ) } \frac{v^p}{d^p_{\phi (\Omega ) }}\left( p V_{\Omega }[\phi,\psi ] -  {\rm div}(\psi\circ\phi^{(-1)}) \right)dy,
\end{eqnarray}
for all $\psi\in C^2(\overline{\Omega} ;{\mathbb{R}}^n)$, where $v$ is a positive minimizer for $H_p(\phi (\Omega ))$
normalized  by   the condition $\| v\|_{L^p(\phi (\Omega) ,  d^{-p}_{\phi (\Omega)})}  =1 $.
\end{thm}
{\em Proof.} Since $H_p(\phi(\Omega ))$ depends continuously on $\phi \in C^2(\overline{\Omega} ;{\mathbb{R}}^n)$,  it easily follows that ${\mathcal{H}}_{\Omega} $ is an open set in  $C^2(\overline{\Omega} ;{\mathbb{R}}^n)$. Now let $\phi \in {\mathcal{H}}_{\Omega}$ be fixed. Applying Theorem~\ref{had} to the open set $\phi (\Omega)$ we obtain that the map $H$ is Gateaux differentiable at $\phi $ and that the Gateaux differential is provided by formula (\ref{frech}). By Theorem~\ref{stabmin} and formula (\ref{frech}) we deduce that the Gateaux differential depends continuously on $\phi$. As is well known, this implies that the map $H$ is also Fr\'{e}chet differentiable. \hfill $\Box$


\begin{rem}
It would be natural to simplify  (\ref{frech}) and write a formula involving  surface integrals. In the case of classical eigenvalue problems  this is usually done by integrating repeatedly by parts.   
For example, if we denote by $\lambda_p(\phi)$ the usual first eigenvalue of the  $p$-Laplacian with Dirichlet boundary conditions  on $\phi (\Omega )$ then one obtains the well-known Hadamard-type formula 
$$
d\lambda_p (\phi )[\psi]= (1-p) \int_{\partial \phi (\Omega)}\left|\frac{\partial w}{\partial \vec{n} }\right|^p  \mu\cdot \vec{n}\, d\sigma ,
$$    
where $w$ is the first normalized eigenfunction, $\vec{n}$ is the unit outer normal to $\partial \phi (\Omega)$ and $\mu :=\psi \circ\phi^{(-1)}$ (cf. \cite{garsab, laplap}).
By applying the same method to (\ref{frech})  and making formal computations one would obtain the meaningless formulas 
\begin{eqnarray}  
\lefteqn{ dH_p(\phi )[\psi]=    
-p H_p(\phi(\Omega ) )\int_{\phi (\Omega)} \frac{v^p}{d^{p+1}_{\phi (\Omega ) }} \nabla d_{\phi (\Omega )}\cdot (\mu \circ \tau_{\phi(\Omega )})\, dy 
} 
\nonumber   \\
& &  +\int_{\partial\phi (\Omega)} |\nabla v|^p\mu\cdot  \vec{n} -p|\nabla v|^{p-2} \sum_{i,j=1}^n \frac{\partial v}{\partial y_i}\frac{\partial v}{\partial y_j}\mu_j  (\vec{n})_i    -H_p(\phi(\Omega ) ) \frac{v^p\mu\cdot  \vec{n}}{d^p_{\phi(\Omega )}} \, d\sigma  \nonumber \\
& &= (1-p) \int_{\partial \phi (\Omega)}\left|\frac{\partial v}{\partial \vec{n} }\right|^p  \mu\cdot \vec{n}\, d\sigma       -H_p(\phi(\Omega ) ) \int_{\partial \phi (\Omega)} \frac{v^p}{d^p_{\phi(\Omega )}}\mu\cdot \vec{n}\, d\sigma  \nonumber \\
& & 
  -p H_p(\phi(\Omega ) )\int_{\phi (\Omega)} \frac{v^p}{d^{p+1}_{\phi (\Omega ) }} \nabla d_{\phi (\Omega )}\cdot (\mu \circ \tau_{\phi(\Omega )})\, dy . \label{meaningless}
\end{eqnarray}

It is clear that the integrals in (\ref{meaningless}) are not well-defined, see  (\ref{decay}). 
In order to bypass this problem one may think of interpreting the above integrals as `principal value integrals' associated with an invading sequence of  open sets relatively compact in $\phi (\Omega)$. However, we prefer not to insist on this.  
\end{rem}

\section{Stability estimates via volume}

Let $\Omega$ be a bounded domain in $\R^n$.  We recall that  $\tau_{\Omega}(x)$ denotes the boundary point nearest to $x\in \Omega$, so that $|x-\tau_{\Omega}(x)|=d_{\Omega}(x)$; we note that $\tau_{\Omega}(x)$ is well defined for almost all $x\in\Omega$. We denote by  $T$ the nonlinear operator given formally by
\begin{equation}
Tw(x)= \int_0^1 |w(x+t(\tau_{\Omega}(x)-x))|dt \; .
\label{t}
\end{equation}
Equivalently,
\begin{equation}
Tw(x)= d_{\Omega}(x)^{-1}\int_{L_x}|w| \, ds \; ,
\label{tt}
\end{equation}
where $L_x$ denotes the line segment with endpoints $x$ and $\tau_{\Omega}(x)$.

\begin{lem}
Let $\Omega$ be a bounded domain of class $C^2$ and let $\delta>0$ be such that $d_{\Omega}$ is $C^2$ in $\Omega_{2\delta}$. 
Then for any $r \in [1, \infty [$ there exists $c>0$ depending only on $n$, $\delta$ and ${\rm Inr}(\Omega)$  such that for all $w\in L^r(\Omega)$ with $w=0$ outside $\Omega_{3\delta /2}$ we have
$$
\int_{\Omega}|Tw|^r dx \le c   \int_{\Omega} \log( {\rm Inr}(\Omega) /d_{\Omega} )|w|^r \, dx .
$$
\label{cl}
\end{lem}
{\em Proof.} Given a continuous function $w$ with  ${\rm supp}\, w\subset \Omega_{3\delta /2}$ we have
\begin{eqnarray*}
\|Tw\|_{L^r(\Omega)}^r &=&\int_{\Omega} \Big( \int_0^1 |w(x+t(\tau_{\Omega}(x)-x))|dt  \Big)^r dx \\
&\leq& \int_{\Omega}  \int_0^1 |w(x+t(\tau_{\Omega}(x)-x))|^rdt  \, dx \\
&=& \int_{\Omega} d_{\Omega}(x)^{-1}\int_{L_x}|w|^r ds \, dx .
\end{eqnarray*}
We define the map
\[
\Psi :\partial\Omega \times ]0, {\rm Inr}(\Omega)[ \longrightarrow \R^n \; , \qquad 
\Psi( \bar{x},t) =\bar{x} -t \vec{n}(\bar{x}) \; ,
\]
where $\vec{n}(\bar{x})$ denotes the unit outer normal at $\bar{x}\in\partial\Omega$. Let
\[
g(\bar{x},t) = \frac{1}{t} \int_{ [ \bar{x} , \bar{x} -t\vec{n}(\bar{x})]}|w|^r ds \, .
\]
By the Area Formula \cite[Section 3.3.2]{ev} we have
\begin{eqnarray*}
\int_{\Omega} d_{\Omega}(x)^{-1}\int_{L_x}|w|^r ds \, dx 
&\leq & \int_{\Psi(\partial\Omega \times ]0, {\rm Inr}(\Omega)[)} \sum_{ (\bar{x},t)\in\Psi^{-1}(x)}
g(\bar{x},t) dx \\
&=& \int_{\partial\Omega}\int_0^{{\rm Inr}(\Omega)} g(\bar{x},t) |\det \nabla\Psi(\bar{x},t)|  dt \, dS(\bar{x}) \\
&\leq& c\int_{\partial\Omega}\int_0^{{\rm Inr}(\Omega)} g(\bar{x},t)   dt \, dS(\bar{x}) \\
&=& c \int_{\partial\Omega}\int_0^{{\rm Inr}(\Omega)}
\frac{1}{t} \int_0^t |w(\bar{x}-s\vec{n}(\bar{x}))|^r ds\,  dt \, dS(\bar{x}) \\
&=& c \int_{\partial\Omega}\int_0^{{\rm Inr}(\Omega)}
\log( {\rm Inr}(\Omega) /s )|w(\bar{x}-s\vec{n}(\bar{x}))|^r \, ds \, dS(\bar{x}) \\
&=& c \int_{\partial\Omega}\int_0^{3\delta /2}
\log( {\rm Inr}(\Omega) /s )|w(\bar{x}-s\vec{n}(\bar{x}))|^r \, ds \, dS(\bar{x})\\
&\leq& c \int_{\Omega_{3\delta /2}}
\log( {\rm Inr}(\Omega) /d_{\Omega} )|w|^r \, dx ,
\end{eqnarray*}
where $c$ is a positive constant depending  only on $n$, $\delta$ and ${\rm Inr}(\Omega)$ (see \cite[\S 14.6]{gitr} for details concerning uniform upper and lower bounds for  $ |\det \nabla\Psi|$). This completes the proof. $\hfill\Box$

We recall (cfr. (\ref{decay})) that if $\Omega$ is a bounded domain of class $C^2$ with $H_p(\Omega)<((p-1)/p)^p$, then
the $L^p$ Hardy quotient has a positive minimizer $u\in W^{1,p}_0(\Omega)$ and there exists a constant $K=K_p(\Omega )>0$ such that
\begin{equation}
\label{sym}
u(x) \leq   K d(x)^{\alpha} \,  , \quad x\in\Omega ;
\end{equation}
here $\alpha>(p-1)/p \, $ denotes the largest solution of equation (\ref{alpha}).

Given a bi-Lipschitz map $\phi :\Omega\to \phi(\Omega)$ we
define for $1\leq r < \infty$ the following measure of vicinity of $\phi$ to the identity map:
\begin{equation}\label{vicinity}
\delta_{r,p}(\phi)= \bigg( \int_{\Omega}\log( {\rm Inr}(\Omega) /d_{\Omega} ) \big(  | \nabla\phi -I |^r + | \nabla\phi -I |^{pr} \big)dx\bigg)^{1/r}.
\end{equation}
Moreover, for any $\gamma<1$ we set
\begin{equation}
I_{\gamma}(\Omega ):=\int_{\Omega }\frac{1}{d_{\Omega }^\gamma}dx.
\end{equation}

\begin{thm}
Let $\Omega$ be a bounded domain in ${\mathbb{R}}^n$ of class $C^2$. Let $\delta>0$ be such that $d_{\Omega}$ is of class $C^2$  on $\Omega_{2\delta}$.
Assume that $H_p(\Omega)<((p-1)/p)^p$ and let $\alpha>(p-1)/p \, $ denote the largest solution of (\ref{alpha}).
Then for any $r> 1/( \alpha p -p +1)$ there exists $c>0$ depending only on  $n$, $\delta$, $p$, $r$, $L$,    $H_p(\Omega)$, 
$K_p(\Omega )$, $I_{pr(1-\alpha )/(r-1)}(\Omega)$ and ${\rm Inr} (\Omega)$  such that
\begin{equation}
 H_p(\phi(\Omega)) \leq  H_p(\Omega) +   c \delta_{r,p}(\phi) .
\label{aaa}
\end{equation}
for all  $\phi\in\bLip_L(\Omega)$ satisfying $\phi=I$ on $\Omega\setminus \Omega_{\delta}$.
\label{thm:aaa}
\end{thm}
{\em Proof.}  Let $u\in W^{1,p}_0(\Omega)$ be  a positive minimizer for $H_p(\Omega)$.
We assume that $\int_{\Omega}u^p/d_{\Omega}^p dx=1$ and set $v=u\circ\phi^{(-1)}$. 
We then have (cf. (\ref{b4}))
\begin{eqnarray}
&& R_{\phi(\Omega)}[v] -R_{\Omega}[u]  \label{rara}  \\
&=& \frac{ \int_{\Omega} \Big( |(\nabla u) (\nabla\phi)^{-1}|^p |\det \nabla\phi | -  |\nabla u|^p \Big) dx  + 
\int_{\Omega}|\nabla u|^p dx \; \Big( 1-  \int_{\Omega}\Frac{u^p|\det \nabla\phi |}{d^p_{\phi(\Omega)}(\phi(x))}  dx \, \Big)}
{\int_{\Omega}\Frac{u^p}{d^p_{\phi(\Omega)}(\phi(x))} |\det \nabla\phi | dx}.
\nonumber
\end{eqnarray}
By (\ref{b6}) we have that
\begin{equation}
\label{b7}
\int_{\Omega}\Frac{u^p}{d^p_{\phi(\Omega)}(\phi(x))} |\det \nabla\phi | dx \geq c^{-1} .
\end{equation}
We also note that for any $x\in\Omega$ we have
\begin{equation}
d_{\phi(\Omega)}\phi(x) \leq d_{\Omega}(x) + F_{\phi }(x)d_{\Omega}(x),
\label{eq:f}
\end{equation}
where
\[
F_{\phi }(x)=d_{\Omega}(x)^{-1} |  (\phi - I)(\tau_{\Omega}(x)) -  (\phi - I)(x)    |.
\]
Hence
\begin{eqnarray}
\lefteqn{\frac{1}{d^p_{\Omega}(x)} -\frac{ |\det \nabla\phi | }{d^p_{\phi(\Omega)}(x)} \le \frac{1}{d^p_{\Omega}(x)} 
\left(1- \frac{ |\det \nabla\phi | }{(1+F_{\phi })^p } \right)}\nonumber \\
& & \qquad \le\frac{1}{d^p_{\Omega}(x)}\left( 1-\frac{1-c |\nabla\phi - I|  }{(1+F_{\phi })^p} \right) \le \frac{c}{d_{\Omega}^p}\big(F_{\phi }+F^p_{\phi } +|\nabla\phi -I|\big)
\label{eq:bbc1a}
\end{eqnarray}
provided $c$ is large enough.
From (\ref{rara}), (\ref{b7}) and (\ref{eq:bbc1a}) we conclude that
\begin{eqnarray}
&& H_p(\phi(\Omega)) - H_p(\Omega) \nonumber \\
&\leq& c \int_{\Omega}|\nabla\phi - I | \cdot  |\nabla u|^pdx +cH_p(\Omega)\int_{\Omega}\frac{u^p}{d_{\Omega}^p}(F_{\phi}+F^p_{\phi }+|\nabla\phi - I|)dx \nonumber \\
&=:& A[u].  \label{dada}
\end{eqnarray}

Now let $r> 1/( \alpha p -p +1)$ be fixed.
By (\ref{marshaf1}) and (\ref{marshaf1grad})  we have $u/d_{\Omega}, \nabla u\in L^{rp/(r-1)}(\Omega)$, hence
\[
\int_{\Omega}|\nabla\phi - I | \cdot  |\nabla u|^pdx \leq \|\nabla\phi -I\|_{L^r(\Omega)}\|\nabla u\|_{L^{\frac{rp}{r-1}}(\Omega)}^p,
\]
and similarly,
\[
\int_{\Omega}|\nabla\phi - I| \frac{u^p}{d_{\Omega}^p}dx \leq 
\|\nabla\phi -I\|_{L^{r}(\Omega)}\| u/d_{\Omega}\|_{L^{\frac{rp}{r-1}}(\Omega)}^p .
\]
Let $\phi_{k}$, $k\in{\mathbb{N}}$, be an approximating sequence of smooth maps obtained by standard mollification 
of $\phi$. Then $\phi_k$ converges to $\phi $ pointwise and in $L^q(\Omega )$ for any $1\le q< \infty$.
Since $\phi_k$ is smooth, we easily see that  $|F_{\phi_k}|\le T|\nabla \phi_k-I| $ where $T$ is the operator defined in (\ref{t}). Hence  using Fatou's  Lemma, Lemma \ref{cl} and observing that ${\rm supp}\, |\nabla \phi_k-I|\subset \Omega_{3\delta /2}$ provided $k$ is sufficiently large, we obtain
\begin{eqnarray}
\| F_{\phi}\|^r_{L^r(\Omega)}&\le& \liminf_{k\to \infty}  \| F_{\phi_k}\|^r_{L^r(\Omega)} 
\leq  c\lim_{k\to \infty } \int_{\Omega}\log( {\rm Inr}(\Omega) /d_{\Omega} ) | \nabla\phi_k -I |^r dx \nonumber \\
&=& c  \int_{\Omega}\log( {\rm Inr}(\Omega) /d_{\Omega} ) | \nabla\phi -I |^r dx  \; ,  
\label{f}
\end{eqnarray}
and similarly
\begin{equation}
\label{ff}
\| F^p_{\phi}\|^r_{L^r(\Omega)}\le c \int_{\Omega}\log( {\rm Inr}(\Omega) /d_{\Omega} ) | \nabla\phi -I |^{pr} dx \; .  
\end{equation}
From (\ref{vicinity}), (\ref{f}) and (\ref{ff}) we obtain 
\begin{eqnarray*}
&& \hspace{-1cm}    \int_{\Omega}\frac{u^p}{d_{\Omega}^p}(F_{\phi }+F^p_{\phi } )dx \\
&\leq& c \bigg(  \int_{\Omega}\log( {\rm Inr}(\Omega) /d_{\Omega} ) \big(  | \nabla\phi -I |^r + | \nabla\phi -I |^{pr} \big)dx \bigg)^{1/r}
 \| u/d_{\Omega} \|_{L^{\frac{rp}{r-1}}(\Omega)}^p  \\
&=& c \delta_{r,p}(\phi)  \|u/d_{\Omega}\|_{L^{\frac{rp}{r-1}}(\Omega)}^p.  
\end{eqnarray*}
Combining the above we obtain
\[
 A[u] \leq  c \delta_{r,p}(\phi) ( \|\nabla u\|_{L^{\frac{rp}{r-1}}(\Omega)}^p +  \|u/d_{\Omega}\|_{L^{\frac{rp}{r-1}}(\Omega)}^p)\le c   \delta_{r,p}(\phi)  \|   d_{\Omega}^{\alpha -1}\|_{L^{\frac{rp}{r-1}}(\Omega)}^p,
\]
which completes the proof. $\hfill\Box$

We say that an open set $\Omega $ in ${\mathbb{R}}^n$ is of class $C^2_M$ for some $M>0$ if   it can be described locally by the subgraphs of functions of class $C^2_M$, i.e.,  
the standard  $C^2$-norms of such functions are bounded by $M$.  
Then we introduce an additional definition. 

{\bf Definition.}
Let $V$ be a bounded open cylinder, i.e., a set which in some coordinate system $(\bar{y},y_n)$  has the form $V=W\times ]a,b[$, for some bounded convex open set $W\subset\R^{n-1}$. Let $M,\rho>0$. We say that a bounded domain $\Omega\subset\R^n$ belongs to $\cC^2_M(V,\rho)$ if
$\Omega$ is of class $C^2_M$  and there exists a function $g\in C^2(\overline{W})$ such that
$a+\rho \leq g\leq b$,   $\|g\|_{C^2(\overline{W})}  \le M$, and
\begin{equation}
 \Omega \cap V=\{ (\bar{y},y_n) \; : \;  \bar{y}\in W \, ,\,  a<y_n<g(\bar{y}) \}.
\label{C11}
\end{equation}
\begin{thm}
Let $\Omega$ be a bounded domain in ${\mathbb{R}}^n$ of class $\cC^2_M(V,\rho)$ 
such that $H_p(\Omega)<((p-1)/p)^p$ and let $\alpha>(p-1)/p \, $ denote the largest solution of (\ref{alpha}).
Then for any $s\in ]0,  \alpha p -p +1[$ there exists $c>0$ depending only on
$p$, $M$, $V$, $\rho$,  $s$, $H_p(\Omega)$, $I_{p(1-\alpha )/(1-s)}(\Omega )$ and ${\rm Inr} (\Omega)$ such that
\begin{equation}
\big|  H_p(\tilde\Omega) -  H_p(\Omega) \big| \leq   c|\tilde\Omega\triangle\Omega |^{s},
\label{fuga}
\end{equation}
for all domains $\tilde \Omega$ of class $\cC^2_M(V,\rho)$ such that $\Omega\setminus V=\tilde\Omega\setminus V$ and 
 $|\tilde\Omega\triangle\Omega |<c^{-1}$.
\label{thm:niki}
\end{thm}
{\em Proof.} We use a construction of \cite{BL1}.
Let $g,\tilde g: W\to\R$ be the $C^2$ functions describing $\Omega\cap V$ and $\tilde\Omega\cap V$
according to the definition above.
Let $\eta =  \frac{\rho}{2 (b-a)}$. We define the function $g_0= \min \{g ,\tilde g\}-\eta |g-\tilde g  |$ and the domain
\begin{equation}
\label{graf2}
\Omega_0:= (\Omega\setminus V) \cup \left\{(\bar y,y_n):\ \bar y\in W,\ a<y_n<g_0(\bar y) \right\}.
\end{equation}
We note that by the choice of $\eta$ we have $g_0>a$.
We next define a map $\phi$ on  $\Omega$ as follows: for $x\in\Omega\setminus V$ we set $\phi(x)=x$, while for $x\in\Omega \cap V$
we use the local coordinates $(\bar y,y_n)$ to define
\begin{equation}
\label{graf3bis}
\phi (\bar y,y_n) = \left\{
\begin{array}{ll}
(\bar y,y_n), & {\rm if}\ (\bar y,y_n)\in  \overline{\Omega}_0 \, , \\
\left(\bar y, \tilde g(\bar y)+a(\bar y)(y_n-g(\bar y)) \right),
 & {\rm if}\ (\bar y,y_n)\in  \Omega \setminus   \overline{\Omega}_0 \, ,
\end{array}
\right.
\end{equation}
where
\[
a(\bar y)=\darr{ \frac{\eta}{\eta+1},}{ \mbox{ if } \tilde g(\bar y)\leq g(\bar y),}{ \frac{\eta+1}{\eta},}{ \mbox{ if } 
\tilde g(\bar y)\geq g(\bar y).}
\]
It can then be seen that $\phi(\Omega)=\tilde\Omega$ and $\phi\in \bLip_{ L}(\Omega)$, where $L$ depends only on $M$ and  $\eta$.
It is immediate that
\begin{equation}
\label{andreas}
\delta_{r_1,p}(\phi) \leq c |\tilde\Omega\triangle\Omega |^{\frac{1}{r}}, \qquad 1\leq r_1<r<\infty
\end{equation}
where $c$ depends only on $M$, $\eta$, $p$, $r_1$, $r$ and $I_{p(1-\alpha )/(1-s)}(\Omega )$.

Since $\Omega$ is of class $C^2_M$ there exists $\delta>0$ depending only on $M$ such that $d_{\Omega}$ is of class $C^2$ in $\Omega_{2\delta}$. In order to apply Theorem \ref{thm:aaa} we need that $\phi= I$ on
$\Omega\setminus\Omega_{\delta}$, and for this it suffices to guarantee that $|g-g_0| <\delta$. This will be the case if we establish that $\|g -\tilde{g}\|_{L^{\infty}(W)}$ is small enough.
Since $g,\tilde{g}$ are of class $C^2_M$, by the Gagliardo-Nirenberg interpolation inequality (see \cite[p.~125]{nire}) the $L^{\infty}$ norm of $g - \tilde{g}$ is estimated from above via its $L^1$ norm, which is precisely $| \Omega\triangle\tilde\Omega |$. Thus if $| \Omega\triangle\tilde\Omega |$ is small enough then $\phi =I$ in $\Omega\setminus\Omega_{\delta}$, hence Theorem \ref{thm:aaa} applies. 
Using an intermediate  $r_1 \in ] 1/(\alpha p -p +1) ,r[$ (say the midpoint of the interval) we obtain that for any $r>1/(\alpha p -p +1)$
\begin{equation}
\label{andreas1}
 H_p(\tilde\Omega ) \leq  H_p(\Omega) +   c |\tilde\Omega\triangle\Omega |^{\frac{1}{r}}.
\end{equation}
In order to prove the reverse inequality   we first note that from (\ref{andreas1}) it follows that $H_p(\tilde\Omega ) < ((p-1)/p)^p$ provided $|\tilde\Omega\triangle\Omega |$  is small enough. Moreover, if $r> 1/( \alpha p -p +1)$ and $| \Omega\triangle\tilde\Omega |$ is small enough then $r> 1/( \tilde \alpha p -p +1)$ where $\tilde \alpha>(p-1)/p \, $ denotes the largest solution of $(p-1){\tilde\alpha}^{p-1}(1-\tilde\alpha)=H_p(\tilde \Omega)$.
Then by using $\tilde \Omega $ as reference domain in the  procedure above we obtain  that 
\begin{equation}
\label{andreas2}
 H_p(\Omega ) \leq  H_p(\tilde \Omega) +   c |\tilde\Omega\triangle\Omega |^{\frac{1}{ r}}
\end{equation}
where $c$ additionally depends on $H_p(\tilde \Omega )$, $I_{pr(1-\tilde\alpha )/(r-1)}(\tilde \Omega)$ and ${\rm Inr}(\tilde \Omega )$.  
Since these last quantities can be controlled  via the corresponding quantities related to $\Omega $ and via $M$ and  $V$,  choosing $r=1/s$ we deduce the validity of ({\ref{fuga}}). 
$\hfill\Box$\\

\begin{rem}
Note that the bigger is  $H_p(\Omega)$,  the smaller is the number  $ \alpha p -p +1$ which defines the range of admissible exponents $s$ in (\ref{fuga}), namely $\alpha p -p +1\to 0$ as $H_p(\Omega)\to  ((p-1)/p)^p$. With regard to this, we observe that the proof of Theorem~\ref{thm:niki} relies on  the existence of a minimizer and we recall that a minimizer does not exist in the limiting case $H_p(\Omega)=  ((p-1)/p)^p$.
\end{rem}

{\bf Acknowledgements} This research was initiated in 2012 when the first author visited the Department of Mathematics of the University of Padova in the frame of the {\it Visiting Scientist Program} of the University of Padova. The first author acknowledges the warm hospitality of the Department of Mathematics, University of Padova, during the above period. This research was also supported by the research project ``Singular perturbation problems for differential operators", Progetto di Ateneo of the University of Padova.

\end{document}